\documentclass[11pt,leqno]{article}

\usepackage{mathptmx}
\usepackage{amsmath,amsthm,amssymb}
\usepackage{enumerate} 
\usepackage{graphicx}

\theoremstyle{plain}
\newtheorem{thm}{Theorem}

\theoremstyle{remark}
\newtheorem{rem}[thm]{Remark}

\theoremstyle{definition}
\newtheorem{ex}[thm]{Example}

\makeatletter

\makeatletter

\makeatletter

\makeatletter

\title{\Large\bf Quasi-Monte Carlo methods for Choquet integrals}
\author{Yumiharu Nakano\footnote{
This study is partially supported by JSPS KAKENHI Grant Number 26800079.
}
\\[1em]
        \small{Graduate School of Innovation Management} \\
        \small{Tokyo Institute of Technology} \\
        \small{2-12-1 W9-117 Ookayama 152-8552, Tokyo, Japan}
}

\date{\today}

\begin{document}

\maketitle

\begin{abstract}
We propose numerical integration methods for Choquet integrals where 
the capacities are given by distortion functions of an underlying probability measure. 
It relies on the explicit representation of the integrals for step functions and 
can be seen as quasi-Monte Carlo methods in this framework.  
We give bounds on the approximation errors in terms of the modulus of continuity of 
the integrand and the star discrepancy.   

\begin{flushleft}
\end{flushleft}
\begin{flushleft}
\end{flushleft}
\end{abstract}

%


In this paper, we are concerned with numerical integration for Choquet integrals 
\begin{equation}
\label{eq:1}
 \int f(U)dc_{\psi}:= \int_0^{\infty}c_{\psi}(f(U)>x)dx 
  + \int_{-\infty}^0(c_{\psi}(f(U)>x)-1)dx   
\end{equation}
for continuous functions $f$ on $[0,1]^d$, 
where $U$ is a $[0,1]^d$-valued and uniformly distributed random variable on an 
atomless probability space $(\Omega,\mathcal{F},\mathbb{P})$, the function 
$\psi:[0,1]\to [0,1]$ is increasing and concave such that 
$\psi(0)=0$, $\psi(1)=1$, and 
$c_{\psi}$ is the submodular set function defined by 
$c_{\psi}(A)=\psi(\mathbb{P}(A))$, $A\in\mathcal{F}$.  
We refer to Denneberg \cite{den:1994} for the theory of Choquet integrals. 

The capacities of the form $c_{\psi}(A)=\psi(\mathbb{P}(A))$ as above 
appear in financial risk management. In particular, 
the case $\psi(t)=\min(t,\lambda)/\lambda$, for some $\lambda\in (0,1)$, 
corresponds to the risk measure known as the average-value-at-risk, 
which is also called as the conditional value-at-risk or the expected shortfall in practice. 
We refer to  Artzner et al.~\cite{art-del-ebe-hea:1999}, 
McNeil et al.~\cite{mcn-fre-emb:2005}, 
and F\"{o}llmer and Schied \cite{fol-sch:2004} for details. 

As for numerical integration, several techniques that are analogous to those for 
the linear integral have been studied in the literature. 
See, e.g., \cite{mcn-fre-emb:2005} for Monte Carlo methods, 
and Nakano \cite{nak:2012} for optimal quantization methods. 
However, to the best of our knowledge, quasi-Monte Carlo methods have not been 
examined to Choquet integrals despite of its popularity. 

To find a suitable quasi-Monte Carlo method for (\ref{eq:1}), 
let $\{u_i\}_{i=1}^n$ is a point set in $[0,1]^d$ and consider 
the simple random variable $U^{(n)}$ defined by 
\begin{equation*}
 U^{(n)}=\sum_{i=1}^nu_i1_{A_i}, 
\end{equation*}
where  $\{A_i\}_{i=1}^n\subset\mathcal{F}$ is a partition of $\Omega$ such that 
$\mathbb{P}(A_i)=1/n$, $i=1,\ldots,n$.  
Note that such $A_i$'s exist since $(\Omega,\mathcal{F},\mathbb{P})$ 
is assumed to be atomless.  
If $\{u_i\}$ is uniformly distributed, then we expect 
$$
 \int f(U)dc_{\psi} \approx \int f(U^{(n)})dc_{\psi}. 
$$

Next, recall that the Choquet integral has the 
{\em comonotonicity\/}, i.e.,  for any random variables $X$ and $Y$ 
that are integrable with respect to $c_{\psi}$, we have 
\begin{equation*}
 \int(X+Y)dc_{\psi}=\int Xdc_{\psi}+\int Ydc_{\psi}
\end{equation*}
whenever   
\begin{equation}
\label{eq:2}
 (X(\omega)-X(\omega^{\prime}))(Y(\omega)-Y(\omega^{\prime}))\ge 0 
\end{equation}
for all $(\omega,\omega^{\prime})\in \Omega\times\Omega$ except for a
set of probability zero. 
Two random variables $X$ and $Y$ are said to be comonotone if 
they satisfy (\ref{eq:2}). 

Now observe that for $A,B\in\mathcal{F}$, 
the two indicator functions
$1_A$ and $1_B$ are comonotone if $A\supset B$. 
Thus, for $A_i\in\mathcal{F}$ and $a_i\in\mathbb{R}$, 
$i=1,\ldots,n$, we have 
\begin{equation}
\label{eq:3}
 \int\sum_{i=1}^na_i1_{A_i}dc_{\psi}
 =\sum_{i=1}^na_i\int 1_{A_i}dc_{\psi} 
\end{equation}
provided that $A_1\subset\cdots\subset A_n$ and $a_i\ge 0$, 
$i=1,\ldots,n$. 

Let $\tau:\{1,\ldots,n\}\to\{1,\ldots,n\}$ be such that 
 $f(u_{\tau(1)})\le\cdots\le f(u_{\tau(n)})$. 
Then we have the representation of $f(U^{(n)})$ given by 
\begin{equation*}
 f(U^{(n)})=\sum_{i=1}^nf(u_{\tau(i)})1_{A_{\tau(i)}}
 =f(u_{\tau(1)})+\sum_{i=2}^n(f(u_{\tau(i)})-f(u_{\tau(i-1)}))
  1_{\cup_{k=i}^nA_{\tau(k)}}. 
\end{equation*}
Applying this representation to (\ref{eq:3}), we obtain 
\begin{equation*}
 \int f(U^n)dc_{\varphi} = f(u_{\tau(1)}) + \sum_{i=2}^n(f(u_{\tau(i)})-f(u_{\tau(i-1)}))
  c_{\psi}\left(\bigcup_{k=i}^nA_{\tau(k)}\right). 
\end{equation*}
Consequently, 
since $c_{\psi}(\cup_{k=i}^nA_{\tau(k)})=\psi(\sum_{k=i}^n\mathbb{P}(A_{\tau(k)}))=\psi((n-i+1)/n)$, 
the quantity 
\begin{equation*}
 I^{(n)}(f) := f(u_{\tau(1)}) + \sum_{i=2}^n(f(u_{\tau(i)})-f(u_{\tau(i-1)}))
  \psi\left(\frac{n-i+1}{n}\right)
\end{equation*}
can be an approximation of $\int f(U)dc_{\psi}$. 

To obtain an error bound, we use the star discrepancy 
$D^*(x_1,\ldots,x_m)$ defined by 
\begin{equation*}
 D^*(x_1,\ldots,x_m)=\sup\left\{\frac{1}{m}\sum_{j=1}^m\left|
  1_B(x_j)-\mathrm{Leb}(B)\right|\: :\:  
  B=\prod_{i=1}^d[0,\xi_i), \; \xi_k\in [0,1],\; k=1,\ldots,d\right\} 
\end{equation*}
for a point set $\{x_i\}_{i=1}^m\subset [0,1]^d$,  
where $\mathrm{Leb}$ stands for the Lebesgue measure on $[0,1]^d$. 
We refer to Niederreiter \cite{nie:1992} for the relation between the discrepancy and 
numerical integration. 
Further, let $\rho(g;t)$ be the modulus of continuity of a function $g$ defined by 
\begin{equation*}
 \rho(g;t)=\sup\{|g(x)-g(y)|: |x-y|\le t, \; x,y\in [0,1]^d\}, \quad t\ge 0, 
\end{equation*}
where $|x|$ denotes the max norm of a vector $x$. 
Also, notice that by the concavity of $\psi$, the limit 
$$
 \psi_+^{\prime}(t):=\lim_{s\searrow t}\frac{\psi(s)-\psi(t)}{s-t}, \quad t\in [0,1), 
$$
exists for any $t\in [0,1)$ and decreasing with respect to $t$. 

Then we have the following: 
\begin{thm}
Under the assumptions and notations above, 
if $\rho(f; D^*(u_1,\ldots,u_n)^{1/d})<1$, we have 
\begin{equation*}
 \left|\int f(U)dc_{\varphi} - I^{(n)}(f)\right|
 \le \left(2\max_{u\in [0,1]^d}|f(u)|+4\right) 
  \psi(\rho(f; D^*(u_1,\ldots,u_n)^{1/d})). 
\end{equation*}
Moreover, if $\psi_+^{\prime}(0)<\infty$, then 
\begin{equation*}
 \left|\int f(U)dc_{\varphi} - I^{(n)}(f)\right| 
 \le 4\psi_+^{\prime}(0)\rho(f; D^*(u_1,\ldots,u_n)^{1/d}). 
\end{equation*}
\end{thm}

\begin{rem}
Since we have assumed that $f$ and $\psi$ are continuous with $\psi(0)=0$, 
the quantity $\psi(\rho(f;D^*(u_1,\ldots,u_n)^{1/d}))$ and so the approximation 
error converge to zero, provided that $D^*(u_1,\ldots,u_n)\to 0$ as $n\to\infty$. 
In particular, if $\psi^{\prime}_+(0)<\infty$, the function $f$ is Lipschitz on $[0,1]^d$, 
and $\{u_i\}_{i=1}^n$ is a low-discrepancy point set, i.e., it satisfies 
$$
 D^*(u_1,\ldots,u_n)\le C\frac{(\log n)^d}{n}
$$
for some positive constant $C$, 
then the theorem implies 
$$
 \left|\int f(U)dc_{\varphi} - I^{(n)}(f)\right| 
 \le 4\psi_+^{\prime}(0)|f|_{Lip}C^{1/d}\frac{\log n}{n^{1/d}}, 
$$
where $|f|_{Lip}$ is the Lipschitz constant of $f$. 
\end{rem}

\begin{rem}
In case $\psi(t)=\min(t,\lambda)/\lambda$, we have  
$\psi_{+}^{\prime}(0)=1/\lambda<\infty$. 
\end{rem}

\begin{rem}
If $d=1$ then the constant $4$ in the statement of the theorem is replaced by $1$. 
This can be verified from the proof below and Theorem 2.10 in \cite{nie:1992}. 
\end{rem}


\begin{proof}[Proof of Theorem 1]
By Lemma 4.63 in \cite{fol-sch:2004} we define the Borel probability measure 
$\mu$ on $[0,1]$ by the identity
\begin{equation*}
 \psi_+^{\prime}(t) = \int_{(t,1]}\frac{1}{s}\mu(ds), \quad t\in (0,1). 
\end{equation*}
Then, from Lemma 4.46 and Theorem 4.64 in \cite{fol-sch:2004} it follows that 
\begin{equation}
\label{eq:4}
 \int Xdc_{\psi}=\int_{[0,1]}\frac{1}{\lambda}\inf_{y\in\mathbb{R}}\left(\mathbb{E}
  (X-y)_++\lambda y\right)\mu(d\lambda) 
\end{equation}
for any bounded random variable $X$, where 
$(x)_+=\max(x,0)$ for $x\in\mathbb{R}$. 
Moreover, for each $\lambda\in (0,1)$, the infimum of the integrand in (\ref{eq:4}) is attained by 
\begin{equation*}
 y_X(\lambda) := \inf\{x\in\mathbb{R}: \mathbb{P}(X>x)\le\lambda\}. 
\end{equation*}

Writing (\ref{eq:4}) with $X=f(U), f(U^{(n)})$, we have for $\varepsilon\in [0,1)$, 
\begin{align*}
 &\int f(U)dc_{\psi} - \int f(U^{(n)})dc_{\psi} \\
 &\le\int_{[0,\varepsilon]}\bigg\{\frac{1}{\lambda}\mathbb{E}[(f(U)-y_{f(U^{(n)})}(\lambda))_+ +  
 \lambda y_{f(U^{(n)})}(\lambda)]  
 -\frac{1}{\lambda}\mathbb{E}[(f(U^{(n)})-y_{f(U^{(n)})}(\lambda))_+ 
 + \lambda y_{f(U^{(n)})}(\lambda)] \bigg\}\mu(d\lambda) \\
 &\quad + \int_{(\varepsilon,1]}\bigg\{\frac{1}{\lambda}\mathbb{E}[(f(U)-y_{f(U^{(n)})}(\lambda))_+ 
 + \lambda y_{f(U^{(n)})}(\lambda)]  
  - \frac{1}{\lambda}\mathbb{E}[(f(U^{(n)})-y_{f(U^{(n)})}(\lambda))_+ 
 + \lambda y_{f(U^{(n)})}(\lambda)] \bigg\}\mu(d\lambda). 
\end{align*}
Since $\inf_{y\in\mathbb{R}}((a-y)_++\lambda y)=\lambda a$, $a\in\mathbb{R}$, 
the first term in the equality just above is at most 
$2\max_{u\in [0,1]^d}|f(u)|\mu([0,\varepsilon])$. 
By Fubini's theorem, the second term is equal to 
\begin{equation}
\label{eq:5}
 \int_{[0,1]^d}F^{(n)}(u)du -\frac{1}{n}\sum_{j=1}^nF^{(n)}(u_i), 
\end{equation}
where 
\begin{equation*}
 F^{(n)}(u)=\int_{(\varepsilon,1]}\frac{1}{\lambda}(f(u)-y_{f(U^{(n)})}(\lambda))_+\mu(d\lambda). 
\end{equation*}
By Theorem 1 in Proinov \cite{pro:1988}, 
the quantity (\ref{eq:5}) is bounded by $4\rho(F^{(n)}; D^*(u_1,\ldots,u_n)^{1/d})$. 
Furthermore, it is straightforward to see that 
$\rho(F^{(n)};t)\le \psi_+^{\prime}(\varepsilon)\rho(f;t)$.  
Summarizing the above arguments, we deduce that 
\begin{equation*}
 \int f(U)dc_{\psi} - \int f(U^{(n)})dc_{\psi}
  \le 2\max_{u\in [0,1]^d}|f(u)|\mu([0,\varepsilon]) 
   + 4\psi_{+}^{\prime}(\varepsilon)\rho(f;D^*(u_1,\ldots,u_n)^{1/d}). 
\end{equation*}
A similar argument shows that 
$\int f(U^{(n)})dc_{\psi} - \int f(U)dc_{\psi}$ is bounded by the right-hand side in 
the inequality just above.
Thus, 
\begin{equation}
\label{eq:6}
 \left|\int f(U)dc_{\psi} - I^{(n)}(f)\right| 
  \le 2\max_{u\in [0,1]^d}|f(u)|\mu([0,\varepsilon]) 
   + 4\psi_{+}^{\prime}(\varepsilon)\rho(f;D^*(u_1,\ldots,u_n)^{1/d}). 
\end{equation}

Now, if $\psi_+^{\prime}(0)<\infty$, then we set $\varepsilon =0$ in (\ref{eq:6}) to 
obtain the second assertion of the theorem. 
Otherwise, we use an argument from the proof of Lemma 4.63 in \cite{fol-sch:2004} to 
obtain $\mu([0,\varepsilon])=\psi(\varepsilon)-\varepsilon\psi_+^{\prime}(\varepsilon)
\ge 0$. Therefore, by the choice $\varepsilon = \rho(f;D^*(u_1,\ldots,u_n)^{1/d})$ 
the right-hand side in (\ref{eq:6}) is estimated as  
\begin{equation*}
  2\max_{u\in [0,1]^d}|f(u)|\psi(\rho(f;D^*(u_1,\ldots,u_n)^{1/d}))  
   + 4\psi(\rho(f;D^*(u_1,\ldots,u_n)^{1/d})). 
\end{equation*}
Thus the first assertion of the theorem follows. 
\end{proof}

\begin{figure}[htbp]
\centering
\includegraphics[width=0.7\columnwidth, bb = 0 0 576 432]{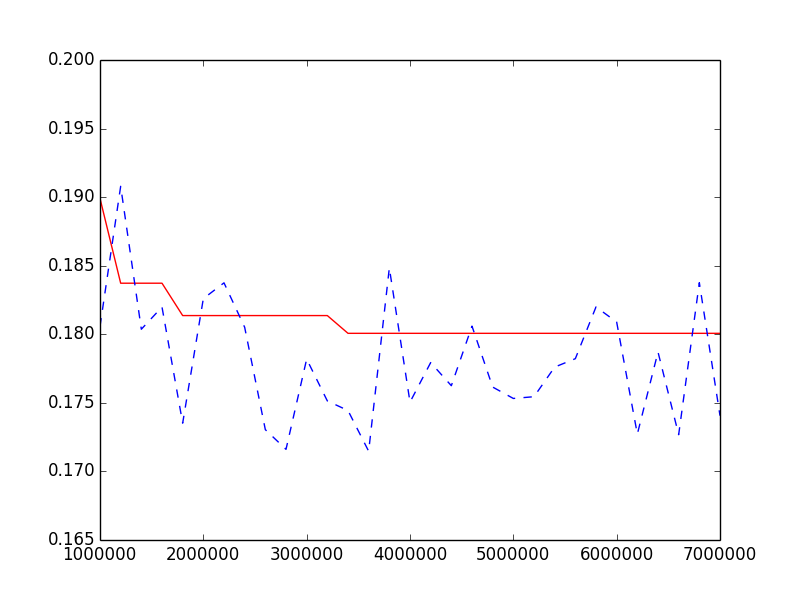}
\caption{Nurimerical integrations of $f(U)$ with quasi-Monte Carlo (solid) 
and Monte Carlo (dashed) methods for $n$ from $10^6$ to $7\times 10^6$.}
\label{fig:1.1}
\end{figure}

\begin{ex}
Here, we present a numerical result in the case of 
$\psi(t)=\min(t,0.05)/0.05$, $t\in [0,1]$, and 
\begin{equation*}
f(u)=\exp\left[-\left\{u^1u^2u^3 + \sin(u^3u^4u^5)\right\}\right], 
\quad u=(u^1,u^2,u^3,u^4,u^5)\in [0,1]^5. 
\end{equation*}
We use the Halton sequence to compute $Q^{(n)}$. 
As a comparison, we take the Monte Carlo method, which is described by 
\begin{equation*}
 I_M^{(n)}(f):= f(U_{\sigma(1)})+\sum_{i=2}^n(f(U_{\sigma(i)})-f(U_{\sigma(i-1)}))
  \psi\left(\frac{n-i+1}{n}\right), 
\end{equation*}
where $\{U_n\}_{n=1}^{\infty}$ is an IID sequence with uniform distribution on $[0,1]^5$  
and $\sigma:\{1,\ldots,n\}\to \{1,\ldots,n\}$ is such that 
$f(U_{\sigma(1)})\le \cdots\le f(U_{\sigma(n)})$. 
Figure \ref{fig:1.1} plots values of $I^{(n)}(f)$ and $I_M^{(n)}(f)$ for 
$n$ from $10^6$ to $7\times 10^6$ with step $2\times 10^5$. 
We can see that $I^{(n)}(f)$ steadily approaches to a true value as $n$ increases, 
whereas the behavior of $I_M^{(n)}(f)$ is still volatile even for 
$n$ larger than $4\times 10^6$.  
\end{ex}

\bibliographystyle{hplain}
\bibliography{../mybib}

\end{document}